\documentclass[10pt]{article}
\usepackage{amssymb}
\usepackage{amsmath,amssymb}\usepackage{cite}\usepackage{mathrsfs}\usepackage[amsmath,thmmarks]{ntheorem}
\usepackage{listings}
\usepackage[titletoc]{appendix}\allowdisplaybreaks

\makeatletter
\renewcommand\theequation{\thesection.\@arabic\c@equation}
\makeatother
\newtheorem{thm}{ Theorem}[section]%
\newtheorem{lem}[thm]{ Lemma}%
\newtheorem{Con}[thm]{Conjecture}%
\newtheorem{Fac}[thm]{Fact}%
\input cyracc.def

\def\f{\noindent}

\def\demo{\f{\bf Proof}\hskip10pt}

\def\qed{\hfill $\Box$}

\begin{document}
\title{\textbf{Two Identities}}
\author{Junyao Pan\,
 \\\\
School of Sciences, University of wuxi, Wuxi, Jiangsu,\\ 214105
 People's Republic of China \\}
\date {} \maketitle

\baselineskip=16pt

\vskip0.5cm

{\bf Abstract:} In the process of studying a conjecture of Holly M. Green and Martin W. Liebeck, we obtain two interesting identities by elementary methods, one is a combinatorial identity, and the other is a number theoretic identity.


\section {Introduction}

Throughout this paper, for all natural number $k,n$ with $0\leq k\leq n$, $P^k_{n}=\frac{n!}{(n-k)!}$ and $\dbinom{n}{k}=\frac{n!}{k!(n-k)!}$, where $0!=1$ and $k!=1\times2\times\cdot\cdot\cdot\times k$ for each positive integer $k$.

It well-known that the combinatorial identity is not only interesting but also important in combinatorial enumerations and probability theory. In deed, many famous combinatorial identities have been found, for example, Vandermonde identity and Li Shanlan identity, as follows:

$$\dbinom{m+n}{k}=\sum\limits_{r=0}^k\dbinom{m}{r}\dbinom{n}{k-r}~{\rm{ and}} ~\dbinom{n+k}{k}^2=\sum\limits_{r=0}^k\dbinom{k}{r}^2\dbinom{n+2k-r}{2k}.$$
So far, many scholars are still studying various combinatorial identities, see \cite{A,CY,M}. However, there are few identities about permutations, and the following is an interesting one:
$$P^k_n-P^k_{n-m}=\sum\limits_{r=1}^k\dbinom{k}{r}P^r_mP^{k-r}_{n-m}.$$

On the other hand, for a finite group $G$, if there exists a positive integer $\lambda$ and two subsets $X,~Y$ of $G$ such that $\lambda G=XY$, that is, there are precisely $\lambda$ pairs $(x,y)\in X\times Y$ such that $g=xy$ for every element $g\in G$, then we say that $X$ and $Y$ divide $G$, and the $Y$ is a code with respect to $X$, in particular, it is a perfect code when $\lambda=1$. These codes have attracted quite a bit of attentions, see for example \cite{E,HXZ,T}. Recently, H. M. Green and M. W. Liebeck gave some relevant results concern the symmetric group $S_n$ and proposed the following conjecture, for details see \cite{GL}.

\begin{Con}\label{pan1-0}\normalfont(\cite[Conjecture~2.3]{GL}~)
Let $n>2k$ and let $j$ be such that $2^j\leq k<2^{j+1}$. Suppose $X=x^{S_n}$ is a conjugacy class in $S_n$. Then $\lambda S_n=XY_k$ if and only if the cycle-type of $x$ has exactly one cycle of length $2^i$ for $0\leq i\leq j$ and all other cycles have length at least $k+1$.
\end{Con}
In the process of investigating Conjecture\ \ref{pan1-0}, we attributed it to some combinatorial identity problems. As we feel that Conjecture\ \ref{pan1-0} is true as well as a lot of calculations, we proposed the following conjecture which is very useful in confirming Conjecture\ \ref{pan1-0}.
\begin{Con}\label{pan1-1}\normalfont
Let $k,m,n$ be three positive integers with $1<2k\leq n-m$ and $1<2k\leq m$. Then for each $1\leq i\leq k$, we have
$$P^{k-i}_k(P^i_{n-k}-P^i_{n-k-m})=m\sum\limits_{r=1}^{i}\sum\limits_{s=0}^{k-i}\dbinom{i}{r}\dbinom{k-i}{s}P^{r-1}_{m-r-s-1}P^s_{r+s-1}P^{i-r}_{n-k+r+s-m}P^{k-i-s}_{k-r-s}.$$
\end{Con}
Unfortunately, we can not give a proof of Conjecture\ \ref{pan1-1}. However, we confirmed the special case for $i=k$ and so we obtained the following elementary while interesting permutation identity which may have not been discovered before.
\begin{thm}\label{pan1-2}\normalfont
Let $m,n$ be two positive integers with $m<n$. Then for each integer $k$ such that $1<2k\leq m$ and $k\leq n-m$, we have $P^k_{n}-P^k_{n-m}=m\sum\limits_{r=1}^{k}\dbinom{k}{r}P^{r-1}_{m-r-1}P^{k-r}_{n-m+r}$.
\end{thm}

Recall that for all positive integer $n$, there exist $a_s=1$ and $a_0,a_1,...,a_{s-1}$ are $0$ or $1$ such that $n=a_{0}1+a_{1}2+a_{2}2^2+\cdot\cdot\cdot+a_{s}2^{s}$, which is called the 2-adic representation of $n$. In fact, the 2-adic representation of numbers was introduced by G. W. Leibnizhas, which has important applications in computing technology. In studying Conjecture\ \ref{pan1-0}, we also got the following interesting property.
\begin{thm}\label{pan1-3}\normalfont
Let $k,\eta,a$ be three positive integers with $1\leq k<\eta=1+2+2^2+\cdot\cdot\cdot+2^j$ and $a=\eta-k$. Suppose that $a_{i_0}1+a_{i_1}2+a_{i_2}2^2+\cdot\cdot\cdot+a_{i_{s_i}}2^{s_i}$ is the 2-adic representation of $k+i$ and $b_i=\sum\limits_{r=t_i}^{s_i}a_{i_r}2^{i_r}$ with $2^{i_{t_i}-1}\leq i<2^{i_{t_i}}$ for each $i=1,2,...,a$. Then $ka=\sum\limits_{i=1}^ab_i$.
\end{thm}

\section {Proofs of Theorems}

In this section, we will use Mathematical Induction to prove the above theorems.

{\bf{Proof of Theorem 1.3:}} Regard the case of $k=1$. Then for all $m\geq2$ and $n-m\geq1$, we see $$P^{1}_{n}-P^{1}_{n-m}=m~{\rm{and}} ~m\sum\limits_{r=1}^{1}\dbinom{1}{r}P^{r-1}_{m-r-1}P^{k-r}_{n-m+r}=m\dbinom{1}{1}P^{0}_{m-2}P^{0}_{n-m+1}=m,$$
and thus Theorem\ \ref{pan1-2} holds for $k=1$. Proof by induction on $k$. So we assume that Theorem\ \ref{pan1-2} is true until $k$ for all $m,n$ with $m\geq2k$ and $n-m\geq k$. 

Next we start to prove the case of $k+1$. In particular, the existence of this situation means that $2(k+1)\leq m$ and $n-m\geq k+1$. For convenience, we set $$\Psi=P^{k+1}_{n}-P^{k+1}_{n-m}~{\rm{and}} ~ \Gamma=m\sum\limits_{r=1}^{k+1}\dbinom{k+1}{r}P^{r-1}_{m-r-1}P^{k+1-r}_{n-m+r}.$$
Noticing that $\Gamma=m\sum\limits_{r=1}^{k}\dbinom{k+1}{r}P^{r-1}_{m-r-1}P^{k+1-r}_{n-m+r}+mP^k_{m-k-2}$. In addition, base on the fact that $$\dbinom{k+1}{r}=\dbinom{k}{r}+\dbinom{k}{r-1}~{\rm{and}} ~P^{k+1-r}_{n-m+r}=(n-m+r)P^{k-r}_{n-m+r-1},$$ it follows that
$$\Gamma=m\sum\limits_{r=1}^{k}[\dbinom{k}{r}+\dbinom{k}{r-1}]P^{r-1}_{m-r-1}(n-m+r)P^{k-r}_{n-m+r-1}+mP^k_{m-k-2}.$$
Thus we have $\Gamma=\Gamma_1+\Gamma_2+mP^k_{m-k-2}$, where $$\Gamma_1=m\sum\limits_{r=1}^{k}\dbinom{k}{r}(n-m+r)P^{r-1}_{m-r-1}P^{k-r}_{n-m+r-1},\Gamma_2=m\sum\limits_{r=1}^{k}\dbinom{k}{r-1}(n-m+r)P^{r-1}_{m-r-1}P^{k-r}_{n-m+r-1}.$$ Moreover, we note that $$\Gamma_1=m(n-m)\sum\limits_{r=1}^{k}\dbinom{k}{r}P^{r-1}_{m-r-1}P^{k-r}_{n-m+r-1}+m\sum\limits_{r=1}^{k}\dbinom{k}{r}rP^{r-1}_{m-r-1}P^{k-r}_{n-m+r-1}.$$
Since $2(k+1)\leq m$ and $n-m\geq k+1$, then by Induction hypothesis we deduce that
$$m(n-m)\sum\limits_{r=1}^{k}\dbinom{k}{r}P^{r-1}_{m-r-1}P^{k-r}_{n-m+r-1}=(n-m)(P^k_{n-1}-P^k_{n-m-1}),$$
and therefore, we arrive at $$\Gamma=(n-m)(P^k_{n-1}-P^k_{n-m-1})+m\sum\limits_{r=1}^{k}\dbinom{k}{r}rP^{r-1}_{m-r-1}P^{k-r}_{n-m+r-1}+\Gamma_2+mP^k_{m-k-2}.$$
On the other hand, as with $P^{k+1}_{n}=nP^k_{n-1}$ and $P^{k+1}_{n-m}=(n-m)P^k_{n-m-1}$, we derive that $$\Psi=nP^k_{n-1}-(n-m)P^k_{n-m-1}=(n-m)(P^k_{n-1}-P^k_{n-m-1})+mP^k_{n-1}.$$

Now we start to confirm $\Psi=\Gamma$. So far, we have seen that $$\Psi=\Gamma~{\rm{ ~if~ and ~only~ if}}~
m\sum\limits_{r=1}^{k}\dbinom{k}{r}rP^{r-1}_{m-r-1}P^{k-r}_{n-m+r-1}+\Gamma_2=mP^k_{n-1}-mP^k_{m-k-2}.$$
And we substitute $\Gamma_2$ in the above equation, it can be concluded that $\Psi=\Gamma$ if and only if
$$\sum\limits_{r=1}^{k}[\dbinom{k}{r}+\dbinom{k}{r-1}]rP^{r-1}_{m-r-1}P^{k-r}_{n-m+r-1}+(n-m)\sum\limits_{r=1}^{k}\dbinom{k}{r-1}P^{r-1}_{m-r-1}P^{k-r}_{n-m+r-1}=P^k_{n-1}-P^k_{m-k-2}.$$
Note that $[\dbinom{k}{r}+\dbinom{k}{r-1}]r=(k+1)\dbinom{k}{r-1}$. Thus we infer that 
$$\Psi=\Gamma~{\rm{ ~if~ and ~only~ if}}~(n-m+k+1)\sum\limits_{r=1}^{k}\dbinom{k}{r-1}P^{r-1}_{m-r-1}P^{k-r}_{n-m+r-1}=P^k_{n-1}-P^k_{m-k-2}.$$
Since $k\leq n-m$ and $2(k+1)\leq m$, we have $2k\leq n-m+k+1$ and $k\leq m-k-2$, and then by Induction hypothesis we see that $$P^k_{n-1}-P^k_{m-k-2}=(n-m+k+1)\sum\limits_{r=1}^{k}\dbinom{k}{r}P^{r-1}_{n-m+k-r}P^{k-r}_{m-k+r-2}.$$
Hence, we deduce that $$\Psi=\Gamma~{\rm{ if~ and ~only ~if}} ~\sum\limits_{r=1}^{k}\dbinom{k}{r-1}P^{r-1}_{m-r-1}P^{k-r}_{n-m+r-1}=\sum\limits_{r=1}^{k}\dbinom{k}{r}P^{r-1}_{n-m+k-r}P^{k-r}_{m-k+r-2}.$$
Indeed, it is clear that $\dbinom{k}{r-1}P^{r-1}_{m-r-1}P^{k-r}_{n-m+r-1}=\dbinom{k}{k-r+1}P^{k-r}_{n-m+r-1}P^{r-1}_{m-r-1}$. Additionally, we observe that
$$\dbinom{k}{k-r+1}P^{k-r}_{n-m+r-1}P^{r-1}_{m-r-1}=\dbinom{k}{k-r+1}P^{(k-r+1)-1}_{n-m+k-(k-r+1)}P^{k-(k-r+1)}_{m-k+(k-r+1)-2},$$ and therefore, we deduce that $$\sum\limits_{r=1}^{k}\dbinom{k}{r-1}P^{r-1}_{m-r-1}P^{k-r}_{n-m+r-1}=\sum\limits_{r=1}^{k}\dbinom{k}{r}P^{r-1}_{n-m+k-r}P^{k-r}_{m-k+r-2}.$$
We have thus proved this theorem by induction.    \qed

{\bf{Proof of Theorem 1.4:}} Consider the case of $j=1$. In this case, $k=1$ or $k=2$. If $k=1$, then $a=2$, $b_1=2$ and $b_2=0$, and thus $ka=b_1+b_2=2$; if $k=2$, then $a=1$ and $b_1=2$, and thus $ka=b_1=2$. So Theorem\ \ref{pan1-3} holds for $j=1$. Assume that Theorem\ \ref{pan1-3} is true until $j=m$ with $m>1$. We claim that Theorem\ \ref{pan1-3} is also true for $j=m+1$. Next we divide into four cases to prove our claim.

Case 1: $1\leq k<1+2+2^2+\cdot\cdot\cdot+2^m$. Pick $b$ such that $k+b=1+2+2^2+\cdot\cdot\cdot+2^m$. Then we see $k+b+2^{m+1}=1+2+2^2+\cdot\cdot\cdot+2^m+2^{m+1}$, in other words, $a=b+2^{m+1}$ and $ka=kb+k2^{m+1}$. By Induction hypothesis, we see $kb=\sum\limits_{i=1}^{b}b_i$. In addition, we note that $b_{b+1}=\cdot\cdot\cdot=b_{b+k}=2^{m+1}$ and $b_{b+k+1}=\cdot\cdot\cdot=b_{b+2^{m+1}}=0$. Thus we have $ka=k(b+2^{m+1})=\sum\limits_{i=1}^{a}b_i$, as desired.

Case 2: $k=1+2+2^2+\cdot\cdot\cdot+2^m$. In this case, we see $a=2^{m+1}$ and $b_i=2^{m+1}$ for $i=1,2,...,2^{m+1}-1$, and $b_{2^{m+1}}=0$. Thus we have $\sum\limits_{i=1}^{a}b_i=(2^{m+1}-1)2^{m+1}=ka$, as desired.

Case 3: $k=2^{m+1}$. In this situation, we have $a=1+2+2^2+\cdot\cdot\cdot+2^m$ and $b_i=2^{m+1}$ for $i=1,2,...,2^{m+1}-1$, and thus $\sum\limits_{i=1}^{a}b_i=(2^{m+1}-1)2^{m+1}=ka$, as desired.

Case 4: $2^{m+1}<k<1+2+2^2+\cdot\cdot\cdot+2^{m+1}$. Since $k+a=1+2+2^2+\cdot\cdot\cdot+2^{m+1}$ and $k>2^{m+1}$, we have $ka=(k-2^{m+1})a+2^{m+1}a$ and $b_i\geq2^{m+1}$. Put $T_i=b_i-2^{m+1}$. Note that $1\leq k-2^{m+1}<1+2+2^2+\cdot\cdot\cdot+2^m$. Then by Case 1, we see $(k-2^{m+1})a=\sum\limits_{i=1}^{a}T_i$. Furthermore, we see $\sum\limits_{i=1}^{a}T_i=\sum\limits_{i=1}^{a}b_i-2^{m+1}a$, and thus $\sum\limits_{i=1}^{a}b_i=\sum\limits_{i=1}^{a}T_i+2^{m+1}a=ka$, as desired.

Hence, the theorem is proven by induction.    \qed

\section{Acknowledgement}

We are very grateful to the anonymous referees for their useful suggestions and comments.

\end{document}